\documentclass[11pt,a4paper]{amsart}

\title[The residue current in codimension three]{The residue current of a codimension three complete intersection}

\author{H\aa kan Samuelsson}\thanks{Author supported by a Post Doctoral Fellowship from the Swedish
Research Council. Part of the work was also done when the author was visiting the ESI in Vienna Oct-Dec 2005,
then supported by a SVeFUM Post Doctoral Fellowship, and he wishes to thank the organizers of the ESI project
'Complex Analysis, Operator Theory and Applications to Mathematical Physics' for their hospitality}
\address{Department of Mathematics, University of Wuppertal, Gaussstrasse 20, 42119 Wuppertal, Germany}
\email{hasam@math.chalmers.se}

\usepackage[english]{babel}
\usepackage{amsmath}
\usepackage{amsthm,amssymb,latexsym}
\usepackage{t1enc}
\usepackage{fancyheadings}
\usepackage{mathrsfs}
\usepackage{graphicx}

\newtheorem{proposition}{Proposition}
\newtheorem{theorem}[proposition]{Theorem}
\newtheorem{lemma}[proposition]{Lemma}

\theoremstyle{definition}

\newtheorem{remark}[proposition]{Remark}

\newcommand{\C}{\mathbb{C}}
\newcommand{\debar}{\bar{\partial}}
\newcommand{\D}{\mathscr{D}}
\newcommand{\R}{\mathbb{R}}
\newcommand{\Z}{\mathcal{Z}}
\newcommand{\V}{\mathcal{V}}
\newcommand{\N}{\mathbb{N}}

\newcommand{\CP}{\mathbb{CP}}
\newcommand{\X}{\mathcal{X}}
\newcommand{\I}{\mathscr{I}}

\newcommand{\real}{\mathfrak{R}\mathfrak{e}}

\def\newop#1{\expandafter\def\csname #1\endcsname{\mathop{\rm #1}\nolimits}}
\newop{span}

\begin{document}
\nocite{*}
\bibliographystyle{plain}

\begin{abstract}
Let $f_1$, $f_2$, and $f_3$ be holomorphic functions on a complex manifold 
and assume that the common zero set of the $f_j$ has maximal codimension.
We prove that the iterated Mellin transform of the residue integral has an analytic
continuation to a neighborhood of the origin in $\C^3$. We prove also that 
the natural regularization of the residue current converges 
unrestrictedly. \\
{\bf MSC}: 32A27; 32C30
\end{abstract}

\maketitle
\thispagestyle{empty}

\section{Introduction}
Let $X$ be an $n$-dimensional complex manifold and $f=(f_1,\ldots,f_m)$ a
holomorphic mapping $X\rightarrow \C^m$ such that the common zero set 
$V_f=\{f=0\}$ has codimension $m$. 
In \cite{CH} Coleff and Herrera were able to associate a certain $(0,n)$-current
to $f$, which has proven to be a good notion of a multi variable residue 
of $f$. They defined their current $R^f$, the Coleff-Herrera 
residue current, as follows. For a test form $\varphi\in\D_{n,n-m}(X)$,
consider the residue integral
\[
I_f^{\varphi}(\epsilon)=
\int_{T_{\epsilon}}\frac{\varphi}{f_1\cdots f_m},
\]
where $T_{\epsilon}$ is the tube 
$\{|f_1|^2=\epsilon_1,\ldots,|f_m|^2=\epsilon_m\}$. Coleff and Herrera
proved that the limit of the residue integral as $\epsilon$ tends to zero along 
a so called admissible path exists and defines the action of a $(0,m)$-current on the 
test form $\varphi$. (The case $m=1$ is due to Herrera and Lieberman, \cite{HL}.)
The limit along an admissible path here means that $\epsilon$ tends to zero 
along a path in the first orthant such that $\epsilon_j/\epsilon_{j+1}^k$ tends
to zero for all $k\in \N$ and $j=1,\ldots,m-1$. The Coleff-Herrera residue current 
has many desirable properties. For instance, it is supported on $V_f$, it has the 
standard extension property, which more or less means that it has no mass concentrated
on the singular parts of $V_f$, and it satisfies the duality property that a 
holomorphic function $h$ on $X$ annihilates it if and only if $h$ 
belongs to the ideal generated by $f$. The duality property is due to 
Dickenstein-Sessa, \cite{DS}, and Passare, \cite{Pdr}, independently.
It is natural to ask if the restriction to limits along admissible paths is
necessary. It actually is and the first example showing this was found by Passare
and Tsikh, \cite{PTmotex}. Björk later realized that this indeed is the typical case, \cite{jeb1};
see also Pavlova, \cite{Pavlova}.
The Coleff-Herrera definition is in this sense quite unstable and one could try to 
look for more stable ones. One step in this direction was taken by Passare in
\cite{Crelle} where he introduced the following regularized version of the residue
integral. Let $\chi_j$ be smooth function on $[0,\infty]$ taking the value $0$
at $0$ and $1$ at $\infty$. (Actually, Passare considered functions identically
$0$ close to $0$ and identically $1$ close to $\infty$.) The regularized residue
integral is then the volume integral
\begin{equation}\label{eqI1}
\int \frac{\debar\chi_1(|f_1|^{2}/\epsilon_1)\wedge\cdots\wedge
\debar\chi_m(|f_m|^{2}/\epsilon_m)}{f_1\cdots f_m}\wedge \varphi.
\end{equation}
Note that in the (not allowed) case when all the $\chi_j$ are the characteristic
function of $[1,\infty]$ we get back the residue integral. It follows from 
Coleff's and Herrera's result that the limit of \eqref{eqI1} along admissible paths 
exists and equals the limit along admissible paths of the residue integral but Passare
proves a more general result. In fact, he proves that for almost all parabolic 
paths, $\epsilon(\delta)=(\delta^{a_1},\ldots,\delta^{a_m})$, the limit of 
\eqref{eqI1} exists. Here "almost all" means that one has to impose finitely many linear
conditions $\sum a_jb_j\neq 0$ in order to asure convergence.
In the special case when $m=2$ the author was able to prove that \eqref{eqI1} actually 
depends Hölder continuously on $\epsilon$ in the closed first quarter, \cite{hasamCR},
\cite{hasamJFA}. I this paper we generalize this result to the codimension three case and 
we prove
\begin{theorem}\label{HS2}
Let $f_1$, $f_2$, and $f_3$ be holomorphic functions on a complex manifold 
$X$ of dimension $n$ and assume that the common zero set of the 
$f_j$ has maximal codimension. Let also $\chi_1$, $\chi_2$, and $\chi_3$
be smooth functions on $[0,\infty]$ taking the value $0$ at $0$ and $1$
at $\infty$ and denote $\chi_j(|f_j|^2/\epsilon_j)$ by $\chi_j^{\epsilon}$.
Then, for any test form $\varphi \in \D_{n,n-2}(X)$, the integral
\[
\int \frac{\chi_1^{\epsilon}\debar \chi_2^{\epsilon}\wedge \debar 
\chi_3^{\epsilon}}{f_1\cdot f_2\cdot f_3}\wedge \varphi
\]
depends Hölder continuously on $\epsilon$ in the closed first octant.
\end{theorem}
The value at the origin is the current $PR^2[1/f]$ (acting on $\varphi$)
in Passare's notation from \cite{Crelle} and it is a $\debar$-potential 
to the Coleff-Herrera residue current. That \eqref{eqI1} converges unrestrictedly
in the case $m=3$ thus follows from Theorem \ref{HS2} by applying it to 
$\debar$-exact test forms.

Another approach to the Coleff-Herrera residue current based on analytic continuation
of currents, a technique with roots in the works of Gelfand and Shilov, \cite{GS}, and 
Atiyah, \cite{At}, has been considered by several authors, e.g., Yger, \cite{Y},
Berenstein, Gay, and Yger, \cite{BGY}, and Passare and Tsikh, \cite{PTcanada},
\cite{PTmotex}. Computing the Mellin transform of the residue integral one obtains
\begin{equation}\label{eqI2}
\int\frac{\debar |f_1|^{2\lambda_1}\wedge\cdots\wedge \debar |f_m|^{2\lambda_m}}{
f_1\cdots f_m}\wedge \varphi,
\end{equation}
where $\lambda_1\ldots,\lambda_m$ are complex parameters with large real parts,
see e.g.\ \cite{BGY} or \cite{PTcanada}. It is proved in \cite{Y} and \cite{BGY} 
that the restriction
of \eqref{eqI2} to any complex line of the form $\lambda(t_1,\ldots,t_m)$, $t_j\in\R_{>0}$,
can be analytically continued to a neighborhood of the origin. Moreover, the
value at the origin equals the Coleff-Herrera residue $R^f.\varphi$. It is also proved in
\cite{BGY} that in the case when $m=2$, \eqref{eqI2} can be analytically continued to a 
neighborhood of the origin as a function of two complex variables. It is generally believed, 
but not yet fully proved, that this holds for arbitrary $m$. Our second main result confirms 
this conjecture for $m=3$.

\begin{theorem}\label{HS1}
Let $f_1$, $f_2$, and $f_3$ be holomorphic functions on a complex manifold 
$X$ of dimension $n$, and assume that the common zero set of the 
$f_j$ has maximal codimension.
Then, for any test form $\varphi \in \D_{n,n-2}(X)$, the holomorphic function
\[
(\lambda_1,\lambda_2,\lambda_3)\mapsto
\int \frac{|f_1|^{2\lambda_1}\debar |f_2|^{2\lambda_2}\wedge \debar 
|f_3|^{2\lambda_3}}{f_1\cdot f_2\cdot f_3}\wedge \varphi,
\]
originally defined when $\real \,\lambda_j$, $j=1,2,3$, is large enough, has a
holomorphic continuation to a neighborhood of $\lambda=0$ in $\C^3$
and the value at the origin equals $PR^2[1/f].\varphi$.
\end{theorem}

We will not be concerned with it in this paper but we
also mention a third, and very successful way to gain stability in the definition 
of the Coleff-Herrera residue current introduced by Passare, Tsikh, and Yger in \cite{PTY}. 
They use the Bochner-Martinelli kernel as a blueprint for the definition of the residue 
instead of the Cauchy kernel. One can view their Bochner-Martinelli type current as the limit 
of a certain average of the residue integral. The advatage is that this averaging process
reduces the number of parameters to just one and it is then proved in \cite{PTY}
that the limit as this single parameter tends to zero exists. It is even 
true when $f$ does not define a complete intersection! However, it is non trivial to
prove that the obtained current actually equals the Coleff-Herrera current.
Quite recently, Andersson put the ideas in \cite{PTY} into an algebraic framework
and introduced more general currents of the Cauchy-Fantappiè-Leray type, which have 
been useful in applications, see e.g.\ \cite{A1}, \cite{A2}, \cite{A3}.

The disposition of the paper is as follows. In the next section we settle the
notations for some frequently appearing objects and we discuss the main elements
of the proofs of Theorems \ref{HS2} and \ref{HS1}. In section \ref{ex} we compute an 
example showing that the codimension three case is different from the codimension two case.
Section \ref{preliminarier} contains some technical results about the normal 
crossings case. At the end of the section we also prove a combinatorial algebra type result
which will enable us to use the assumption that $f$ defines a complete intersection
efficiently. In the last section, Section \ref{bevis}, we prove our main theorems.
However, we only prove Theorem \ref{HS1} in detail since the proofs are almost
identical.
 
\section{Notations and an overview of the proof}
We will have to use Hironaka's theorem, \cite{H}, to resolve singularities locally. 
It gives us for any sufficiently small open set $U\subset X$ a complex manifold $\X$ and 
a proper holomorphic map $\pi\colon \X\rightarrow U$ with the properties that 
$\Z_f:=\{\pi^*f_1\cdot \pi^*f_2\cdot \pi^*f_3=0\}$ has normal crossings and $\pi$ restricted
to $\X\setminus \Z_f$ is a biholomorphism. The varieties we will be most interested in are 
the varieties $V_f:=\{f_1=f_2=f_3=0\}$ and $Z_f:=\{f_1f_2f_3=0\}$ in $X$, and their 
total transforms, $\V_f:=\{\pi^*f_1=\pi^*f_2=\pi^*f_3=0\}$ and 
$\Z_f=\{\pi^*f_1\cdot \pi^*f_2\cdot \pi^*f_3=0\}$. Varieties in calligraphic letters are 
always varieties in the resolution manifold $\X$. Moreover, varieties denoted by $Z$ (or
$\Z$) are always varieties of codimension $1$ and a holomorphic function as a subscript
means the zero variety of that holomorphic function, e.g., $Z_{f_1}:=\{f_1=0\}$
($\Z_{f_1}:=\{\pi^*f_1=0\}$). Varieties of higher codimensions are denoted by
$V$ (or $\V$). Occasionally we will encounter varieties for which this nomenclature is 
not efficiently applicable and we will then use more ad hoc notations. 
That the variety $\Z_f$ has normal crossings in $\X$ means that locally on $\X$
one can find holomorphic coordinates $z$ such that $\pi^*f_j=z^{a_j}\tilde{f}_j$,
$j=1,2,3$, where the $\tilde{f}_j$ are non vanishing holomorphic functions. We will 
call $z_k$ a simple factor if $z_k$ divides precisely one of the monomials
$z^{a_j}$. 
The following analytic sheaves on $\X$ will be referred to frequently; the sheaves of
holomorphic $k$-forms, $\Omega^k$, and the subsheaves of them of 
holomorphic $k$-forms vanishing on a normal crossings divisor $\Z$, $\I^k_{\Z}$.
A holomorphic $k$-form, $\alpha$, vanishes on a normal crossings divisor if the
pullback of $\alpha$ (under the inclusion map) to any irreducible component of $\Z$
vanishes. If $z$ are local coordinates such that $\Z$ is the zero set of a monomial,
$z^{a}$, then $\alpha$ vanishes on $\Z$ if and only if $(dz_j/z_j)\wedge \alpha$
is a holomorphic $k+1$-form for any $z_j$ dividing the monomial $z^a$. This, in turn,
holds if and only if, for all $r\geq 1$, 
$(dz_{j_1}/z_{j_1})\wedge\cdots \wedge (dz_{j_r}/z_{j_r}) \wedge\alpha$ is a 
holomorphic $k+r$-form for $z_{j_l}$ dividing $z^a$.

Now, some comments to the proofs of Theorems \ref{HS2} and \ref{HS1}.
After a partition of unity we may assume that
our test form $\varphi$ has support in a neighborhood $U$ such that it exists 
a Hironaka resolution of singularities $\pi\colon \X\rightarrow U$ as described above.
We then pull our integral 
back to the resolution manifold $\X$ and as explained above, we find local holomorphic 
coordinates $z$ on $\X$ such that 
$\pi^*f_j=z^{a_j}\tilde{f}_j$, $j=1,2,3$. After a partition of unity on $\X$
one is then able to start computing. From a computational point of 
view it is of course easier if
one could arrange so that the $\tilde{f}_j\equiv 1$. This is possible if the 
integer vectors $a_j$ are linearly independent, see e.g.\ \cite{Crelle}. If one 
restricts to limits along admissible paths, as in \cite{CH}, or allowed parabolic paths,
as in \cite{Crelle}, one will encounter only charts on $\X$ where the $a_j$ are linearly
independent. However, in the general case one will encounter also charts where the 
$a_j$ are linearly dependent, so called charts of resonance. This is precisely what happens
in the Passare-Tsikh example, \cite{PTmotex}. Charts of resonance can therefore be seen 
as the reason for the discontinuity of the residue integral. 
In codimension two the author showed in \cite{hasamCR} and \cite{hasamJFA} how the 
charts of resonance can be handled when one considers the regularized residue integral
\eqref{eqI1}. The main tool is Proposition $11$ in \cite{hasamJFA} and it can be generalized. 
The general version is Proposition \ref{tekniskprop} below but we have
omitted the proof since it is a  straightforward generalization of the proof of 
Proposition $11$ in \cite{hasamJFA}. The presence of charts of resonance is therefore not
a problem when we prove Theorems \ref{HS2} and \ref{HS1}. The main difficulty is another
problem that is harder to handle when $m\geq 3$ then when $m=2$. 
On the resolution manifold $\X$ the functions
$\pi^*f_j$ almost never define a complete intersection and for arbitrary test forms in 
$\D_{n,n-m}(\X)$ the corresponding residue integral will in general be discontinuous. 
One has to use that the test forms we consider are of the special form
$\pi^*\varphi$ because, in such test forms, the information that $f$ defines a complete intersection 
in $X$ is somehow coded. When $m=2$ it is quite easy to extract this information.
Actually, it follows from a degree argument, see e.g.\ \cite{CH}, 
that $\pi^*\varphi$ vanishes on all components of codimension $1$ of the variety
$\V_f$, i.e., on the exceptional divisor. This vanishing is then
seen to be enough to get the results in codimension two. In codimension three however,
$\pi^*\varphi$ will in general not vanish on the exceptional divisor with the consequence
that it is not a local problem on the resolution manifold $\X$ to prove Theorems 
\ref{HS2} and \ref{HS1}. In the next section we give a simple example showing this.
So one has to work a little more to extract the information hidden
in $\pi^*\varphi$ when $m=3$.
The degree argument is still very useful though. With the aid of the slightly technical
Lemma \ref{CHlemma} it enables us to locally modify the
test form $\pi^*\varphi$, without affecting the integral, so that the modified test 
form has good enough vanishing properties on the exceptional divisor. The process,
however, produces also a global term which requires some additional attention.

\section{An example}\label{ex}

We consider an example showing that proving analyticity of the Mellin transform and
continuity of the regularized residue integral are not local problems on the 
resolution manifold. We will look at the integral
\begin{equation}\label{eq0}
\int \frac{|x_1|^{2\lambda_1}\debar |x_2|^{2\lambda_2}\wedge \debar
|x_3|^{2\lambda_3}}{x_1x_2x_3}\wedge \varphi(x) dx\wedge d\bar{x}_1
\end{equation}
in $\C^3$, where $\varphi$ is a function defined as follows. Let 
$\phi$, $\varphi_2$ and $\varphi_3$ be smooth functions on 
$\C$ with support close to the origin but non vanishing there, and put
$\varphi_1=\partial\phi/\partial\bar{z}$. We define $\varphi(x)$ to be the function
$\varphi_1(x_1)\varphi_2(x_2)\varphi_3(x_3)$ in $\C^3$. First of all, note that \eqref{eq0}
is the Mellin transform of a residue integral by the choice of $\varphi$;  
we can move the $\debar_1$ in front of $\phi$ to $|x_1|^{2\lambda_1}$ by an integration by parts.
Secondly, \eqref{eq0} equals 
\[
\int \frac{|x_1|^{2\lambda_1} |x_2|^{2\lambda_2}|x_3|^{2\lambda_3}}{x_1x_2x_3} 
\varphi_1 \frac{\partial \varphi_2}{\partial \bar{x}_2}
\frac{\partial \varphi_3}{\partial \bar{x}_3}dx\wedge d\bar{x}
\]
after two integrations by parts, from which we see that \eqref{eq0} is analytic at $\lambda=0$.
Now we blow up $\C^3$ along the $x_1$-axis and study the pullback of \eqref{eq0} to this 
manifold. Let $\pi\colon \C\times \mathcal{B}_0\C^2\rightarrow \C^3$ be the blow up map.
In the natural coordinates $z$ and $\zeta$ on $\C\times \mathcal{B}_0\C^2$ it looks like
\[
\pi(z_1,z_2,z_3)=(z_1,z_2,z_2z_3),
\]
\[
\pi(\zeta_1,\zeta_2,\zeta_3)=(\zeta_1,\zeta_2\zeta_3, \zeta_2).
\]
Since $\varphi$ has support close to the origin, $\pi^*\varphi$ has support close to
$\pi^{-1}(0)=\{z_1=z_2=0\}\cup\{\zeta_1=\zeta_2=0\}\cong \CP^1$. Note that $z_3$ and 
$\zeta_3$ are natural coordinates on this $\CP^1$ and choose a partition of unity,
$\{\rho_1,\rho_2\}$ on $\mbox{supp} (\pi^*\varphi)$ such that $\mbox{supp}
(\rho_1)\subset \{|z_3|<2\}$ and $\mbox{supp} (\rho_2)\subset \{|\zeta_3|<2\}$.
The pullback of \eqref{eq0} under $\pi$ now equals
\[
\int \frac{|z_1|^{2\lambda_1}\debar |z_2|^{2\lambda_2}\wedge \debar |z_2z_3|^{2\lambda_3}}{z_1z_2^2z_3}
\wedge \rho_1(z)\varphi_1(z_1)\varphi_2(z_2)\varphi_3(z_2z_3)z_2dz\wedge d\bar{z}_1
\]
\[
-\int \frac{|\zeta_1|^{2\lambda_1}\debar |\zeta_2\zeta_3|^{2\lambda_2}\wedge
\debar |\zeta_2|^{2\lambda_3}}{\zeta_1\zeta_2^2\zeta_3}\wedge
\rho_2(\zeta)\varphi_1(\zeta_1)\varphi_2(\zeta_2\zeta_3)\varphi_3(\zeta_2)\zeta_2
d\zeta\wedge d\bar{\zeta}_1.
\]
We know that this sum (difference) is analytic at $\lambda=0$ but we will check that non of the terms 
are. We consider the first term. It is easily verified that it can be written as
\[
\frac{\lambda_2}{\lambda_2+\lambda_3}
\int \frac{|z_1|^{2\lambda_1}\debar |z_2|^{2(\lambda_2+\lambda_3)}\wedge 
\debar |z_3|^{2\lambda_3}}{z_1z_2z_3}
\wedge \rho_1(z)\varphi_1(z_1)\varphi_2(z_2)\varphi_3(z_2z_3)dz\wedge d\bar{z}_1.
\]
We denote this integral, with the coefficient $\lambda_2/(\lambda_2+\lambda_3)$ removed,
by $I(\lambda)$. After two integrations by parts one sees that $I(\lambda)$ is analytic at the
origin, and so $\lambda_2I(\lambda)/(\lambda_2+\lambda_3)$ is analytic at the origin if and 
only if $I(\lambda)$ vanishes on the hyperplane $\lambda_2+\lambda_3=0$. In particular
we must have that $I(0)=0$. But $I(0)$ can be computed using Cauchy's formula, and one obtains
$I(0)=-(2\pi i)^3\phi(0)\varphi_2(0)\varphi_3(0)\neq 0$. Hence, proving analyticity of the 
Mellin transform of the residue integral is not a local problem on the blown up manifold.
The same example can also be used to see that continuity of the regularized residue integral
is not a local property on the resolution manifold.

\begin{remark}
This example could be a little confusing. The variable $z_1$ just appears as a "dummy variable" 
in the computations above, to which nothing interesting happens. This indicates that it is not
a local problem on the resolution manifold to prove
analytic continuation of \eqref{eqI2} (or an unconditional limit of \eqref{eqI1}) already in 
the case $m=2$. Actually, if one pulls \eqref{eqI2} back to a resolution manifold and then
starts computing one will encounter global problems already for $m=2$. But
analytic continuation of 
\begin{equation}\label{eqI22}
(\lambda_1,\lambda_2)\mapsto \int 
\frac{|f_1|^{2\lambda_1}\debar |f_2|^{2\lambda_2}}{f_1\cdot f_2}\wedge \varphi, \,\,\,\,
\varphi \in \D_{n,n-1}(X),
\end{equation} 
implies analytic continuation of \eqref{eqI2} for $m=2$, and proving analytic continuation
of \eqref{eqI22} is a local problem on the resolution manifold. Thus, the analyticity
problem can always be reduced to a local problem on the resolution manifold when $m=2$ but, 
and this is the point, for $m\geq 3$ this is not always possible. 
\end{remark}

\section{Preliminary lemmas}\label{preliminarier}
The first proposition in this section is a  straightforward generalization of 
Proposition $11$ in \cite{hasamJFA} so we omit the proof. 

\begin{proposition}\label{tekniskprop}
Let $\psi_j$, $j=1,\ldots,m$, be strictly positive smooth functions on
an open set $\Omega \subset \C^n$ and let $a_j=(a_{j1},\ldots,a_{jn})$, 
$j=1,\ldots,m$, be multiindices. Let also $\chi_j\in C^{\infty}([0,\infty])$ 
be zero at zero. Then, for any test form $\phi \in \D_{n,n}(\Omega)$, the integral
\[
\int \frac{\chi_1(\psi_1 |z^{a_1}|^2/\epsilon_1)\cdots
\chi_1(\psi_m |z^{a_m}|^2/\epsilon_m)}{z^{a_1}\cdots z^{a_m}}
\wedge \phi
\]
depends Hölder continuously on $\epsilon=(\epsilon_1,\ldots,\epsilon_m)$ in the closed
first orthant.
\end{proposition}

The following two lemmas more or less reduce the proofs of Theorems
\ref{HS2} and \ref{HS1} to a study of the pullback of the test form to 
the resolution manifold with the result that the two theorems can be treated
almost identically.
 
\begin{lemma}\label{tillrlemma1}
Let $\tilde{f}_j$, $j=1,\ldots,m$, be non vanishing holomorphic functions on
an open set $\Omega \subset \C^n$ and let $a_j=(a_{j1},\ldots,a_{jn})$, 
$j=1,\ldots,m$, be multiindices. Assume that the test form
$\varphi\in \D_{n,n-r}(\Omega)$ has the property that 
$(d\bar{z}_k/\bar{z}_k)\wedge \varphi \in \D_{n,n-r+1}(\Omega)$ for all 
non simple factors $z_k$ dividing some monomial $z^{a_j}$ with 
$1\leq j \leq r$. Then the integral 
\[
\int \frac{\debar |z^{a_1}\tilde{f}_1|^{2\lambda_1}\wedge \cdots \wedge
\debar |z^{a_r}\tilde{f}_r|^{2\lambda_r}\cdot |z^{a_{r+1}}\tilde{f}_{r+1}|^{2\lambda_{r+1}}\cdots
|z^{a_m}\tilde{f}_m|^{2\lambda_m}}{z^{a_1}\tilde{f}_1\cdots z^{a_m}\tilde{f}_m}\wedge \varphi
\]
has an analytic continuation to neighborhood of $\lambda=0$ in $\C^m$. 
\end{lemma}

\begin{lemma}\label{tillrlemma2}
Let $\tilde{f}_j$ and $a_j$, $j=1,\ldots,m$, and $\varphi$ be as in Lemma 
\ref{tillrlemma1} and let $\chi_j\in C^{\infty}([0,\infty])$ be zero at zero.
Then the integral
\[
\int \frac{\debar \chi^{\epsilon}_1\wedge \cdots \wedge 
\debar \chi^{\epsilon}_r \cdot\chi^{\epsilon}_{r+1}\cdots 
\chi^{\epsilon}_m}{z^{a_1}\tilde{f}_1\cdots z^{a_m}\tilde{f}_m}\wedge \varphi,
\]
where $\chi^{\epsilon}_j=\chi_j(|z^{a_j}\tilde{f}_j|^2/\epsilon_j)$, depends Hölder 
continuously on $\epsilon=(\epsilon_1,\ldots,\epsilon_m)$ in the closed first orthant. 
\end{lemma}

\begin{proof}[Proof of Lemmas \ref{tillrlemma1} and \ref{tillrlemma2}.]
It is well known that integrals of the form 
\begin{equation}\label{eq1}
\int \frac{|z^{a_1}\tilde{f}_1|^{2\lambda_1}\cdots |z^{a_m}\tilde{f}_m|^{2\lambda_m}}{
z^{a_1}\tilde{f}_1\cdots z^{a_m}\tilde{f}_m}\wedge \phi
\end{equation}
have an analytic continuation to a neighborhood of $\lambda=0$ in $\C^m$
without any assumptions on the $(n,n)$-test form $\phi$, see e.g.\ \cite{A1}.
Moreover, by Proposition \ref{tekniskprop} we have Hölder continuity in the first 
orthant for integrals like
\begin{equation}\label{eq2}
\int \frac{\chi^{\epsilon}_1\cdots \chi^{\epsilon}_m}{
z^{a_1}\tilde{f}_1\cdots z^{a_m}\tilde{f}_m}\wedge \phi
\end{equation}
for all $(n,n)$-test forms $\phi$. Using the assumption on our test form $\varphi$
we will reduce the computations of the integrals in the Lemmas \ref{tillrlemma1} 
and \ref{tillrlemma2} to sums of integrals of the forms \eqref{eq1} and \eqref{eq2}
respectively. This is done in more or less the same way in both cases. We start by writing
every $\debar$ as the sum $\debar = \debar_1+\cdots+\debar_n$ and splitting up the integrals 
into sums accordingly. We deal with an expression $\debar_k |z^{a_j}\tilde{f}_j|^{2\lambda_j}$,
respectively $\debar_k \chi_j(|z^{a_j}\tilde{f}_j|^2/\epsilon_j)$, as follows. 
If $z_k$ is a non simple factor dividing the monomial $z^{a_j}$ we let $\debar_k$ act,
obtaining
\[
\lambda_j|z_{a_j}\tilde{f}_j|^{2\lambda_j}\big(a_{jk}\frac{d\bar{z}_k}{\bar{z}_k}+
\overline{\frac{\partial_k \tilde{f}_j}{\tilde{f}_j}}\big), 
\]
respectively
\[
\tilde{\chi}_j(|z^{a_j}\tilde{f}_j|^2/\epsilon_j)\big(a_{jk}\frac{d\bar{z}_k}{\bar{z}_k}+
\overline{\frac{\partial_k \tilde{f}_j}{\tilde{f}_j}}\big),
\]
where $\tilde{\chi}_j(t)=t\chi'_j(t)$. Note that $\tilde{\chi}_j$ is zero at zero and 
smooth on $[0,\infty]$, since $\chi'_j(t)\in \mathcal{O}(1/t^2)$ as $t\rightarrow \infty$,
and hence satisfies the properties required by Proposition \ref{tekniskprop}. 
The assumption on our test form $\varphi$
means that, for any expression $d\bar{z}_k/\bar{z}_k$ arising in this way,
$(d\bar{z}_k/\bar{z}_k)\wedge \varphi$ is again a test form. More generally, it is easy to see,
e.g.\ by making a Taylor expansion á la Lemma $6$ in \cite{hasamJFA} of the cofficients of $\varphi$,
that the assumption on $\varphi$ implies that 
\[
\frac{d\bar{z}_{k_1}}{\bar{z}_{k_1}}\wedge \cdots \wedge
\frac{d\bar{z}_{k_p}}{\bar{z}_{k_p}}\wedge \varphi
\]
is a test form if the $z_{k_l}$ are non simple factors each dividing some monomial
$z^{a_j}$ with $1\leq j \leq r$. Hence, all singular forms 
$d\bar{z}_k/\bar{z}_k$ arising from the action of $\debar_k$, where $k$ is such that
$z_k$ a non
simple factor, can be incorporatad in the test form. On the other hand, if
$z_k$ is a simple factor dividing the monomial $z^{a_j}$ we do not let 
$\debar_k$ act on $|z^{a_j}\tilde{f}_j|^{2\lambda_j}$, respectively
$\chi_j(|z^{a_j}\tilde{f}_j|^2/\epsilon_j)$. Instead we then integrate
this $\debar_k$ by parts. Since $z_k$ is a simple factor it does not divide any
monomial other then $z^{a_j}$ and so, after the integration by parts, $\debar_k$ 
will not encounter any monomial containing $z_k$ as a factor and hence not produce 
the singular expression $d\bar{z}_k/\bar{z}_k$. Hence, the integrals in Lemmas
\ref{tillrlemma1} and \ref{tillrlemma2} can be written as sums of integrals of the 
form \eqref{eq1} and \eqref{eq2} respectively, concluding the proof.
\end{proof}

The rest of this section is devoted to a proof of the following lemma. It
will enable us to use the fact that we have a complete intersection on the
original manifold in an efficient way when we do computations on the blown up one,
where we in general do not have complete intersection.

\begin{lemma}\label{CHlemma}
Consider the monomials $\sigma=z_1^{a_1}\cdots z_{r-1}^{a_{r-1}}$ and 
$\tau=z_r\cdots z_s$ and let $\alpha$ be a holomorphic $k$-form such that
$d\sigma \wedge \alpha \in \I^{k+1}_{Z_\tau}$. Then there exists a holomorphic
$k$-form $\alpha'$ such that
\begin{itemize}
\item[(i)] $d\sigma \wedge \alpha'=0$,
\item[(ii)] $\alpha'\in \I^{k}_{Z_\sigma}$, and
\item[(iii)] $\alpha-\alpha'\in \I^k_{Z_\tau}$.
\end{itemize}
\end{lemma}

The lemma will follow from the next one, which says that property (i) implies
property (ii), and Proposition \ref{taylorprop}, which should be compared to 
Lemma $6$ in \cite{hasamJFA}.

\begin{lemma}\label{lemma1}
Let $\sigma$ be a monomial and $\alpha$ a holomorphic $k$-form. 
Then $(d\sigma/\sigma)\wedge \alpha \in \Omega^{k+1}$ if and only if
$\alpha \in \I^k_{Z_\sigma}$. 
\end{lemma}

\begin{proof}
We have by definition that $\alpha \in \I^k_{Z_\sigma}$ if and only if 
$(dz_j/z_j)\wedge \alpha \in \Omega^{k+1}$ for all $z_j$ dividing $\sigma$
and the "if"-part of the lemma is clear. For the "only if"-part we will use induction
on the number of coordinate functions $z_j$ dividing $\sigma$. (One could also
use Proposition \ref{taylorprop} but we choose to give a direct argument.) If just one 
coordinate function divides $\sigma$ then we are done, again by definition.
We therefore assume that we have proved the "only if"-direction for 
$p-1$ coordinate functions dividing $\sigma$. Now let 
$z^a=z_1^{a_1}\cdots z_{p-1}^{a_{p-1}}$ and assume that 
$(d(z^az_p^{a_p})/z^az_p^{a_p})\wedge \alpha \in \Omega^{k+1}$. 
We then write $\alpha=\alpha'+dz_p\wedge\alpha''$,
where $\alpha'$ and $\alpha''$ do not contain any $dz_p$. Then
\[
\frac{d(z^az_p^{a_p})}{z^az_p^{a_p}}\wedge \alpha=
\frac{d(z^a)}{z^a}\wedge \alpha'+
dz_p\wedge (\frac{\alpha'}{z_p}-\frac{d(z^a)}{z^a}\wedge \alpha'')\in \Omega^{k+1}.
\]
Since the first term on the right hand side does not contain any $dz_p$ it 
follows that both terms on the right hand side are in $\Omega^{k+1}$. Then by the
induction hypothesis, $\alpha'\in \I^k_{z^a=0}$. Moreover, 
$\alpha'/z_p-(d(z^a)/z^a)\wedge \alpha''\in \Omega^k$ since $\alpha'$ and $\alpha''$
do not contain any $dz_p$. But then $\alpha'$ must be divisible with $z_p$ and
$(d(z^a)/z^a)\wedge \alpha''$ must be in $\Omega^k$ (since $\alpha'/z_p$ is smooth 
in $z_1,\ldots,z_{p-1}$ and $(d(z^a)/z^a)\wedge \alpha''$ is smooth in $z_p$).
Thus, $\alpha'\in \I^k_{z^az_p=0}$,
and again by the induction hypothesis, $\alpha''\in \I^k_{z^a=0}$. Hence,
$\alpha=\alpha'+dz_p\wedge \alpha''\in \I^k_{z^az_p=0}$, finishing the induction step.
\end{proof}

\begin{proposition}\label{taylorprop}
Consider the monomial $\tau=z_r\cdots z_s$ and the corresponding variety 
$Z_{\tau}$ in $\C^n$.
Denote the index set $\{r,\ldots,s\}$ by $I$ and let $I(j)$ denote an arbitrary
subset of $I$ with precisely $j$ elements fewer than $I$. Let also $V_{I(j)}$ and
$Z_{I(j)}$ be 
the varieties $\cap_{i\in I(j)}\{z_i=0\}$ and
$\cup_{i\in I\setminus I(j)}\{z_i=0\}$ respectively. (In this notation
$Z_{\tau}=Z_{I(s-r+1)}$.) For any holomorphic $k$-form 
$\omega$ we then let $\omega_{I(j)}$ denote the holomorphic $k$-form obtained from 
$\omega$ by first pulling $\omega$ back to $V_{I(j)}$ and then extending
constantly to $\C^n$. 
Now, let $\alpha$ be a holomorphic $k$-form and put $\alpha^1=\alpha-\alpha_I$ 
and recursively, $\alpha^{i+1}=\alpha^i-\sum_{I(i)}\alpha^i_{I(i)}$. Then
\begin{equation}\label{tayloreq}
\alpha=\alpha_I+\sum_{I(1)}\alpha^1_{I(1)}+\cdots +
\sum_{I(s-r)}\alpha^{s-r}_{I(s-r)}+\alpha^{s-r+1},
\end{equation}
where $\alpha^i_{I(i)}\in \I^k_{Z_{I(i)}}$ and $\alpha^{s-r+1}\in \I^k_{Z_{\tau}}$.
\end{proposition}

\begin{proof}
Using induction over the number of coordinate functions dividing the monomial
$\tau$ it is easy to see that \eqref{tayloreq} holds and so, what remains is to see that 
the $\alpha^i_{I(i)}$ have the correct vanishing properties. We fix $r$ and $s$
with $r\leq s$ and we show that 
$\alpha^i_{I(i)}\in \I^k_{Z_{I(i)}}$ for $i=1,\ldots,s-r+1$, again with induction. 
Note that $Z_{\tau}=Z_{I(s-r+1)}$ and that $\alpha^{s-r+1}=\alpha^{s-r+1}_{I(s-r+1)}$.
First we put $i=1$. We have $\alpha^1=\alpha-\alpha_I$ and so 
$\alpha^1_I=\alpha_I-\alpha_I=0$. Now, if $I(1)=I\setminus \{j\}$, then the pullback 
of $\alpha^1_{I(1)}$ to $\{z_j=0\}$ equals $\alpha^1_I=0$. Hence,
$\alpha^1_{I(1)}\in \I^k_{Z_{I(1)}}$. 
For the induction step, assume that $\alpha^{p-1}_{I(p-1)}\in \I^k_{Z_{I(p-1)}}$.
We have $\alpha^p=\alpha^{p-1}-\sum_{I(p-1)}\alpha^{p-1}_{I(p-1)}$ by definition.
If $I'$ is a fixed set of the type $I(p-1)$ we get that
\[
\alpha^p_{I'}=\alpha^{p-1}_{I'}-\sum_{I(p-1)}(\alpha^{p-1}_{I(p-1)})_{I'}=
\alpha^{p-1}_{I'}-\alpha^{p-1}_{I'}=0.
\]
The second equality follows from the induction hypothesis since if $I(p-1)\neq I'$
then $I'$ contains at least one index $j$ not in $I(p-1)$. Then, since 
$\{z_j=0\}\subset Z_{I(p-1)}$, the induction hypothesis implies that
$\alpha^{p-1}_{I(p-1)}\in \I^k_{z_j=0}$, which in turn gives
that $(\alpha^{p-1}_{I(p-1)})_{I'}=0$.
Now, let $I''$ be a set of the type $I(p)$. We can write $I''=I'\setminus \{j\}$ 
for non unique $I'$ of the type $I(p-1)$ and $j$. Then the pullback of 
$\alpha^p_{I''}$ to $\{z_j=0\}$ equals $\alpha^p_{I'}=0$. Repeating this
for all possible decompositions $I''=I'\setminus \{j\}$ of $I''$ we get 
$\alpha^p_{I''}\in \I^k_{Z_{I''}}$ finising the induction step.
\end{proof}

\begin{proof}[Proof of Lemma \ref{CHlemma}.]
Property (ii) follows from property (i) according to Lemma \ref{lemma1}. We claim that 
\[
\alpha'=\alpha_I+\sum_{I(1)}\alpha^1_{I(1)}+\cdots +
\sum_{I(s-r)}\alpha^{s-r}_{I(s-r)},
\]
where we have used the notations from Proposition \ref{taylorprop}, 
has the properties (i) and (iii). That $\alpha'$ has the property (iii)
is part of the statement of Proposition \ref{taylorprop} so we only
need to check that it has property (i). 
By assumption we have that
\begin{equation}\label{eq3}
d\sigma\wedge \alpha'=
d\sigma\wedge\alpha_I+d\sigma\wedge\sum_{I(1)}\alpha^1_{I(1)}+\cdots +
d\sigma\wedge\sum_{I(s-r)}\alpha^{s-r}_{I(s-r)}\in \I^{k+1}_{Z_\tau}.
\end{equation}
If we pull back $d\sigma\wedge\alpha'$ to $V_I=\{z_r=\cdots =z_s=0\}$ we get
by Proposition \ref{taylorprop} that the pullback of 
$d\sigma\wedge\alpha_I$ is zero since all other terms on the right hand side
of \eqref{eq3} vanish on this set by Proposition \ref{taylorprop}. 
But $d\sigma\wedge\alpha_I$ is independent
of all $z_j$ and $dz_j$ with $j=r,\ldots,s$ and so $d\sigma\wedge\alpha_I=0$
in $\C^n$. Next we pull $d\sigma\wedge\alpha'$ back to 
$V_{I\setminus \{r\}}=\{z_{r+1}=\cdots =z_s=0\}$. Then, using that $d\sigma\wedge\alpha_I=0$
and Proposition \ref{taylorprop}, we get that the pullback  
of $d\sigma\wedge \alpha^1_{I\setminus \{r\}}$ to this set is zero. But 
$d\sigma\wedge \alpha^1_{I\setminus \{r\}}$ is independent of all $z_j$ and 
$dz_j$ with $j=r+1,\ldots,s$ and thus vanishes in all of $\C^n$. Continuing in this
way, running through the indices, then pulling back to varieties of dimension
$+1$ and running through
pairs of indices and so on, we eventually obtain that 
$d\sigma\wedge \alpha'=0$ in $\C^n$.
\end{proof}

\section{Proof(s) of Theorems \ref{HS2} and \ref{HS1}}\label{bevis}

We are now in a position to prove Theorems \ref{HS2} and \ref{HS1}. Apart from
using Lemma \ref{tillrlemma1} when proving Theorem \ref{HS1} and Lemma
\ref{tillrlemma2} when proving Theorem \ref{HS2} the proofs are almost
identical and we choose to focus on Theorem \ref{HS1}. Any differences will
be pointed out explicitly.

\begin{proof}[Proof of Theorems \ref{HS2} and \ref{HS1}.]
After a preliminary partition of unity on $X$ we may assume that the test form has as 
small support as we want. Moreover, from \cite{BGY} we 
are done if the support of the test form does not intersect $V_f=\{f_1=f_2=f_3=0\}$. 
(In the case of Theorem \ref{HS2} this follows from \cite{hasamJFA}.) Assume 
therefore that $\varphi$ has support in a small neighborhood $O$ of 
a point $x\in V_f$. We may 
assume that $\varphi$ has the form $\varphi=\tilde{\phi}\wedge \overline{\phi}$, where 
$\tilde{\phi}$ is a smooth $(n,0)$-form with support close to $x$ and $\phi$ is 
a holomorphic $(n-2)$-form.
Hironaka's theorem implies that there is a complex $n$-dimensional manifold $\X$ and a 
proper holomorphic map $\pi\colon \X\rightarrow O$ such that the variety 
$\Z_{f}=\{\pi^*f_1\cdot\pi^*f_2\cdot\pi^*f_3=0\}$ has normal crossings in $\X$ and 
$\pi$ is biholomorphic outside $\Z_{f}$. Since $\pi$ is proper and 
$\varphi$ has support close to 
$x\in V_f$, the pullback, $\pi^*\varphi$, has compact support close to 
$\pi^{-1}(x)\subset \V_{f}:=\{\pi^*f_1=\pi^*f_2=\pi^*f_3=0\}$. We also introduce 
the notation $\Z_{12}$ for the variety consisting of the components of codimension
one of $\{\pi^*f_1=\pi^*f_2=0\}$ on which $\pi^*f_3$ does not vanish identically. The 
varieties $\Z_{23}$ and $\Z_{13}$ are defined analogously.
Now, consider a point $p\in \V_f$. Since $\Z_{f}$ has normal crossings we may choose local
holomorphic coordinates $z$ close to $p$ such that $z(p)=0$, 
$\pi^*f_1=z^a=z_1^{a_1}\cdots z_{r-1}^{a_{r-1}}$ and $\pi^*f_j$, $j=2,3$, are 
monomials times non vanishing holomorphic functions.  
Generically, $p$ does not lie on $\Z_{23}$, but if it does, then, 
after possibly renumbering the 
coordinates $z_r,\ldots,z_n$, we can write $\Z_{23}=\{z_r\cdots z_s=0\}$. We can symbolically
let $s<r$ denote the case that $p$ does not lie on 
$\Z_{23}$. Since $f_1,f_2,f_3$ is a regular sequence, $df_1\wedge \phi$ vanishes on
$\{f_2=f_3=0\}$ for degree reasons, and so $d\pi^*f_1\wedge \pi^*\phi=
d(z^a)\wedge \pi^*\phi\in \I^{n-1}_{\Z_{23}}$. By Lemma \ref{CHlemma} we may
therefore choose a holomorphic $n-2$-form, $\alpha_1$, in a neighborhood $U$ of $p$ such
that $d\pi^*f_1\wedge \alpha_1=0$, $\alpha_1\in \I^{n-2}_{\Z_{f_1}}(U)$, and 
$\phi-\alpha_1\in \I^{n-2}_{\Z_{23}}(U)$. 
In case $p$ does not lie on $\Z_{23}$
we can take $\alpha_1=0$ close to $p$. In the same way, perhaps after shrinking
$U$, we can 
find $\alpha_2$ and $\alpha_3$ such that $d\pi^*f_2\wedge \alpha_2=0$,
$\alpha_2\in \I^{n-2}_{\Z_{f_2}}(U)$, and $\pi^*\phi-\alpha_2\in \I^{n-2}_{\Z_{13}}(U)$
and similarily for $\alpha_3$. In this way we get an open covering of $\V_f$ and we choose 
a locally finite subcovering, $\{U_j\}_1^{\infty}$. In each $U_j$ we have holomorphic
$n-2$-forms, $\alpha_1^j$, $\alpha_2^j$, and $\alpha_3^j$ with the properties described 
above. We may assume that $\pi^*\varphi$ has support in
$\cup U_j$ and then, since $\pi^*\varphi$ has compact support, it suffices to take 
finitely many $U_j$ to cover it. For convenience we denote this finite family 
by $\{U_j\}_1^{q}$. Subordinate to this family we choose a partition of unity,
$\{\rho_j\}_1^q$, such that $\sum_1^q\rho_j=1$ on the support of $\pi^*\varphi$.
The integral we are interested in can now be written
\begin{equation}\label{eq5}
\sum_1^q\int \frac{|\pi^*f_1|^{2\lambda_1}\debar |\pi^*f_2|^{2\lambda_2}\wedge \debar 
|\pi^*f_3|^{2\lambda_3}}{\pi^*f_1\cdot \pi^*f_2\cdot \pi^*f_3}\wedge \rho_j\pi^*
\tilde{\phi}\wedge \pi^*\bar{\phi}.
\end{equation}
Since $d\pi^*f_i\wedge \alpha_i^j=0$ in $U_j$ we may replace $\pi^*\bar{\phi}$ with
$\pi^*\bar{\phi}-\bar{\alpha}_2^j-\bar{\alpha}_3^j$ in \eqref{eq5} without affecting 
the integral. After this is done we integrate the $\debar$ in front of 
$|\pi^*f_2|^{2\lambda_2}$ by parts, obtaining
\begin{equation}\label{eq6}
\sum_1^q\int \frac{|\pi^*f_1|^{2\lambda_1} |\pi^*f_2|^{2\lambda_2} \debar 
|\pi^*f_3|^{2\lambda_3}}{\pi^*f_1\cdot \pi^*f_2\cdot \pi^*f_3}\wedge \rho_j
\debar(\pi^*\tilde{\phi}\wedge (\pi^*\bar{\phi}-\bar{\alpha}_2^j-\bar{\alpha}_3^j))
\end{equation}
\begin{equation}\label{eq7}
-\sum_1^q\int \frac{\debar |\pi^*f_1|^{2\lambda_1}|\pi^*f_2|^{2\lambda_2}\wedge \debar 
|\pi^*f_3|^{2\lambda_3}}{\pi^*f_1\cdot \pi^*f_2\cdot \pi^*f_3}\wedge \rho_j\pi^*
\tilde{\phi}\wedge (\pi^*\bar{\phi}-\bar{\alpha}_2^j-\bar{\alpha}_3^j)
\end{equation}
\begin{equation}\label{eq8}
+\int \sum_1^q \frac{|\pi^*f_1|^{2\lambda_1}|\pi^*f_2|^{2\lambda_2} \debar 
|\pi^*f_3|^{2\lambda_3}}{\pi^*f_1\cdot \pi^*f_2\cdot \pi^*f_3}\wedge \debar\rho_j\wedge
\pi^*\tilde{\phi}\wedge (\pi^*\bar{\phi}-\bar{\alpha}_2^j-\bar{\alpha}_3^j).
\end{equation}
We first show that each term in the sum \eqref{eq6} has an analytic continuation 
to a neighborhood of $\lambda=0$. To this end, we fix a $j$ and write 
$\pi^*f_i=z^{a_i}\tilde{f}_i$, $i=1,2,3$, as monomials times non vanishing functions
in $U_j$. (The multiindices $a_i$ and the functions $\tilde{f}_i$ of course
also depend on $j$ but we suppress this dependence to avoid to many subscripts.)
By Lemma \ref{tillrlemma1} it is sufficient to show that 
\begin{equation}\label{eq9}
\frac{d\bar{z}_k}{\bar{z}_k}\wedge 
\debar(\pi^*\tilde{\phi}\wedge (\pi^*\bar{\phi}-\bar{\alpha}_2^j-\bar{\alpha}_3^j))
\end{equation}
is a smooth form for all non simple factors $z_k$ dividing $z^{a_3}$. Assume first that 
$z_k$ also divides $z^{a_2}$. Then we write \eqref{eq9} as
\begin{equation}\label{eq10}
\frac{d\bar{z}_k}{\bar{z}_k}\wedge \pi^*\debar\varphi-
\frac{d\bar{z}_k}{\bar{z}_k}\wedge \debar(\pi^*\tilde{\phi}\wedge
(\bar{\alpha}_2^j+\bar{\alpha}_3^j)).
\end{equation}
Since $\alpha_2^j$ and $\alpha_3^j$ vanish on $\{z_k=0\}$ the last term is a smooth form.
On the other hand, $\debar\varphi$ can be written as a sum of products of 
smooth $(n,0)$-forms and anti-holomorphic $(0,n-1)$-forms. These anti-holomorphic forms 
vanish on $\{f_2=f_3=0\}$ for degree reasons, and hence, their pullback under $\pi$
vanish on $\{z_k=0\}$. Thus, the first term in \eqref{eq10} is also a smooth form.
It remains to check the case when $z_k$ divides $z^{a_1}$ and $z^{a_3}$ but not
$z^{a_2}$, i.e., that $\{z_k=0\}$ is (part of) an irreducible component of
$\Z_{13}$. We then write \eqref{eq9} as
\[
\frac{d\bar{z}_k}{\bar{z}_k}\wedge \debar
(\pi^*\tilde{\phi}\wedge (\pi^*\bar{\phi}-\alpha_2^j))-
\frac{d\bar{z}_k}{\bar{z}_k}\wedge \debar
(\pi^*\tilde{\phi}\wedge \alpha_3^j),
\]
which is a smooth form 
since both $\alpha_3^j$ and $\pi^*\bar{\phi}-\alpha_2^j$ vanish on 
$\Z_{13}$.

Next we show that each term in the sum \eqref{eq7} has an analytic continuation to a 
neighborhood of $\lambda=0$. Now we have a $\debar$ in front of $|\pi^*f_1|^{2\lambda_1}$
and so we may subtract also $\alpha_1^j$ in the test form. By Lemma \ref{tillrlemma1}
it will thus be sufficient to show that 
\begin{equation}\label{eq11}
\frac{d\bar{z}_k}{\bar{z}_k}\wedge \pi^*\tilde{\phi}\wedge (\pi^*\bar{\phi}-
\bar{\alpha}_1^j-\bar{\alpha}_2^j-\bar{\alpha}_3^j)
\end{equation}
is a smooth form for any non simple factor $z_k$ dividing at least one of 
$z^{a_1}$ and $z^{a_3}$. If $z_k$ divides all three monomials $z^{a_1}$,
$z^{a_2}$, and $z^{a_3}$ then $\alpha_1^j$, $\alpha_2^j$, and $\alpha_3^j$
all vanish on $\{z_k=0\}$. But then $\pi^*\phi$ does also, because the 
holomorphic $n-2$-form $\phi$ vanishes on $V_f$ for degree reasons
and hence $\pi^*\phi$ vanishes on $\V_f\supset \{z_k=0\}$. In this case, \eqref{eq11}
is therefore a smooth form. If instead $z_k$ divides precisely two of the monomials,
say for simplicity $z^{a_1}$ and $z^{a_2}$, then $\alpha_1^j$, $\alpha_2^j$, and
$\pi^*\phi-\alpha_3^j$ vanish on $\{z_k=0\}$ and we see again that \eqref{eq11}
is a smooth form.

Finally we prove that \eqref{eq8} has an analytic continuation to 
a neighborhood of $\lambda=0$. If one could choose the $\alpha_i^j$ to agree
on overlaps, i.e., to be independent of $j$, then this would be easy
because in this case,
the only thing depending on $j$ would be $\rho_j$, and since the $\rho_j$ sum up 
to $1$, the $\debar \rho_j$ sum up to $0$, and \eqref{eq8} would be identically $0$.
Maybe it is possible to choose the $\alpha_i^j$ in such a good way but even
if the differences $\alpha_i^j-\alpha_i^k$ are not $0$ on overlaps
they have sufficiently good properties for our purposes. Actually,
$(dz_l/z_l)\wedge (\alpha_i^j-\alpha_i^k)$ is a holomorphic form,
where it is defined, for any non simple factor $z_l$. We verify this for $i=2$.
If $z_l$ divides $\pi^*f_2$ it is clear. On the other hand, if $z_l$ does not
divide $\pi^*f_2$, but
divides $\pi^*f_1$ and $\pi^*f_3$, i.e., is (part of) a component of $\Z_{13}$,
then it follows from the fact that 
$\alpha_2^j-\alpha_2^k=\pi^*\phi-\alpha_2^k-(\pi^*\phi-\alpha_2^j)\in \I^{n-2}_{\Z_{13}}$. 
We will show that \eqref{eq8} is analytic by showing that each point 
in $\mbox{supp} (\pi^*\varphi)$ has a neighborhood
such that if $\rho$ is a smooth compactly supported function in this neighborhood
then the integral of $\rho$ multiplied with the sum in \eqref{eq8}
has an analytic continuation to a 
neighborhood of $\lambda=0$. It will then follow, after another partition of 
unity, that \eqref{eq8} also has. 
Now, consider a point $p\in\mbox{supp} (\pi^*\varphi)$ and let 
$U_{j_1},\ldots,U_{j_r}$ be those sets from our cover which contains $p$.
Let $U$ be a neighborhood of $p$ contained in $\cap_iU_{j_i}$ and such
that the $\rho_j$ sum up to $1$ in $U$ 
and let $\rho$ be a smooth function with support in $U$.
The differences $\beta_1^{kl}=\alpha_1^{j_k}-\alpha_1^{j_l}$ and 
$\beta_2^{kl}=\alpha_2^{j_k}-\alpha_2^{j_l}$ with $1\leq k,l\leq r$ are 
all defined in $U$. If we multiply the sum in \eqref{eq8} by $\rho$ and integrate
we obtain
\begin{equation*}
\int \sum_{k=1}^r \frac{|\pi^*f_1|^{2\lambda_1}|\pi^*f_2|^{2\lambda_2} \debar 
|\pi^*f_3|^{2\lambda_3}}{\pi^*f_1\cdot \pi^*f_2\cdot \pi^*f_3}\wedge 
\rho\debar\rho_{j_k}\wedge
\pi^*\tilde{\phi}\wedge (\pi^*\bar{\phi}-\bar{\alpha}_2^{j_k}-\bar{\alpha}_3^{j_k})=
\end{equation*}
\begin{equation}\label{eq12}
\int \sum_{k=1}^r \frac{|\pi^*f_1|^{2\lambda_1}|\pi^*f_2|^{2\lambda_2} \debar 
|\pi^*f_3|^{2\lambda_3}}{\pi^*f_1\cdot \pi^*f_2\cdot \pi^*f_3}\wedge 
\rho\debar\rho_{j_k}\wedge
\pi^*\tilde{\phi}\wedge (\pi^*\bar{\phi}-\bar{\alpha}_2^{j_1}-\bar{\alpha}_3^{j_1})
\end{equation}
\begin{equation}\label{eq13}
- \sum_{k=1}^r \int\frac{|\pi^*f_1|^{2\lambda_1}|\pi^*f_2|^{2\lambda_2} \debar 
|\pi^*f_3|^{2\lambda_3}}{\pi^*f_1\cdot \pi^*f_2\cdot \pi^*f_3}\wedge 
\rho\debar\rho_{j_k}\wedge
\pi^*\tilde{\phi}\wedge (\bar{\beta}_2^{1k}+\bar{\beta}_3^{1k}).
\end{equation}
The only thing in the sum in \eqref{eq12} that depends on $k$ is $\rho_{j_k}$ and 
since the $\rho_{j_k}$ sum up to $1$ in $U$ the $\debar\rho_{j_k}$ sum up
to $0$, and so \eqref{eq12} is identically $0$. But $(dz_l/z_l)\wedge\beta_i^{1k}$
is a holomorphic form for any non simple factor $z_l$ by the discussion above,
and so by Lemma \ref{tillrlemma1}, each term in the sum \eqref{eq13} has an 
analytic continuation to a neighborhood of $\lambda=0$. The proof is complete.
\end{proof}

\end{document}